# Numerical solution for fractional model of telegraph equation by using q-HATM

## P. Veeresha, D. G. Prakasha[*]

*Department of Mathematics, Karnatak University, Dharwad - 580 003, India.*
*E-mail: viru0913@gmail.com, dgprakasha@kud.ac.in, prakashadg@gmail.com*

**Abstract**

The pivotal aim of the present work is to demonstrate an efficient analytical technique, called q-homotopy analysis transform method (q-HATM) in order to analyse a fractional model of telegraph equations. Numerical examples are illustrated to examine the efficiency of the proposed technique. The numerical solutions are obtained in the form of a series solution. The proposed method manipulates and controls the series solution, which rapidly converges to the exact solution in a short admissible domain in an efficient manner.

## 1. Introduction

In last few years, fractional calculus has become the center of many studies because of its broad applications in mathematics, physics, and engineering such as electromagnetism, fluid mechanics, signal processing, hydrology, fractional kinetics, electrochemistry, viscoelasticity, optics, robotics, biomedical and so on. The use of fractional differentiation for the mathematical modelling of physical phenomena has been widespread in various fields of science and engineering. The fundamental results related with solution of fractional differential equations may be found in books [1-5]. The theory of fractional differential equations helps to translate the reality of nature in a better and systematic manner. The fractional derivative considered is in the Caputo sense. Fractional partial differential equations have been carried out by many authors through several techniques such as iterative Laplace transform method; Adomian decomposition method; homotopy analysis method; and fractional reduced differential transform method.

On the other hand, communication systems play a vital role in many of the real world problems. A typical engineering problem involves the transmission of a signal from one point to another. A transmission media is part of the circuit and represents a physical system that directly propagates the signal between two or more points. Certainly, all the transmission media have loss in signal. In order to optimize the transmission media it is needed to determine signal losses. To evaluate these losses, it is necessary to formulate some kind of equations which can calculate these losses efficiently. Telegraph equations arise in the study of propagation of electrical signals in a cable of transmission line and wave phenomena [6]. The telegraph equations are due to Oliver Heaviside [7] who developed the transmission line model. This model demonstrates that the electromagnetic waves can be reflected on the wire, and that appear wave patterns along the transmission line. The telegraph equations are in terms of voltage ($v$) and current ($i$) for a section of a transmission media and that are applicable in several fields such as wave propagation [8], random walk theory [9], signal analysis [10] and etc. In most of the practical situation, these telegraph equations occur in fractional order, not always in integer order. There have been many analytic and numerical methods are available in the literature to solve fractional order telegraph equations. For details we refer the readers to the papers [11-27] and the references therein. Here, in this paper we consider multi-dimensional telegraph equations of fractional order. The one-dimensional ($1D$) space-time fractional telegraph equation is given by

$$\frac{\partial^{2\alpha}v}{\partial t^{2\alpha}} + 2a\frac{\partial^{\alpha}v}{\partial t^{\alpha}} + b^2 v = \frac{\partial^{2\beta}v}{\partial x^{2\beta}} + f(x,t), \quad 0 < \alpha, \beta \leq 1, \tag{1}$$

subjected to the initial and boundary conditions

$$\left.\begin{array}{l} v(x,0) = \phi_1(x), v_t(x,0) = \phi_2(x), \\ v(0,t) = \psi_1(t), v_x(0,t) = \psi_2(t) \end{array}\right\}.$$

Also, the two-dimensional ($2D$) time-fractional telegraph equation is defined as

$$\frac{\partial^{2\alpha}v}{\partial t^{2\alpha}} + 2a\frac{\partial^{\alpha}v}{\partial t^{\alpha}} + b^2 v = \frac{\partial^2 v}{\partial x^2} + \frac{\partial^2 v}{\partial y^2} + f(x,y,t), \quad 0 < \alpha \leq 1, \tag{2}$$

subjected to the initial and boundary conditions

$$v(x,y,0) = \xi_1(x,y), v_t(x,y,0) = \xi_2(x,y).$$

Similarly, the three-dimensional ($3D$) time-fractional telegraph equation can be given

$$\frac{\partial^{2\alpha}v}{\partial t^{2\alpha}} + 2a\frac{\partial^{\alpha}v}{\partial t^{\alpha}} + b^2 v = \frac{\partial^2 v}{\partial x^2} + \frac{\partial^2 v}{\partial y^2} + \frac{\partial^2 v}{\partial z^2} + f(x,y,z,t), \quad 0 < \alpha \leq 1, \tag{3}$$

subjected to the initial and boundary conditions

$$v(x,y,z,0) = \chi_1(x,y,z), v_t(x,y,z,0) = \chi_2(x,y,z).$$

In Eqs. (1), (2) and (3), $a$ and $b$ denote positive constants. Also, $f, \phi_1, \phi_2, \psi_1, \psi_2, \xi_1, \xi_2, \chi_1$ and $\chi_2$ are known continuous functions.

The aim of this paper is to propose a new analytic technique for fractional space-time telegraph equations. This technique is a combined form of the q-homotopy analysis method (q-HAM) and Laplace transform. This technique is called the q-homotopy analysis transform method (q-HATM) and yields the series solution in a large admissible domain which helps us to control the convergence region of the series solution. The accuracy and efficiency of the proposed method are demonstrated by the five test examples. This method is easy to implement for the multi-dimensional space and time fractional order physical problems emerging in various fields of engineering and science. This method was also successfully applied to fractional model of Fokker-Planks equation by Prakash and Kaur [28], fractional model of vibration equation by Srivastava et al. [29], fractional Fitzhugh – Nagumo equation by Kumar et al. [30], fractional coupled Burger's equation by Singh et al. [31] and many others.

## 2. Preliminaries

In this section, we recall some basic definitions and properties of fractional calculus which will be used in the subsequent sections.

**Definition 1.** The left sided Riemann-Liouville fractional integral of order $\alpha > 0$ of a function $f \in C_\mu$, $\mu \geq -1$ is defined as:

$$J^\alpha f(t) = \frac{1}{\Gamma(\alpha)} \int_0^t \frac{f(\tau) d\tau}{(t-\tau)^{1-\alpha}},$$

$$J^0 f(t) = f(t).$$

**Definition 2.** The left sided Caputo fractional derivative of $f$, $f \in C_{-1}^m$, $m \in \mathbb{N} \cup \{0\}$

$$D_t^\alpha f(t) = \begin{cases} \frac{d^n f(t)}{dt^n}, & \alpha = n \in \mathbb{N}, \\ \frac{1}{\Gamma(n-\alpha)} \int_0^t (t-\tau)^{n-\alpha-1} f^n(\tau) d\tau, & n-1 < \alpha < n, n \in \mathbb{N}. \end{cases}$$

For instance $f(t) = t^\beta$, we have;

$$D_t^\alpha t^\beta = \begin{cases} \dfrac{\Gamma(\beta+1)}{\Gamma(\beta-\alpha+1)} t^{\beta-\alpha}, & n-1 < \alpha \leq n, \beta > n-1, \beta \in \mathbb{R}, \\ 0, & n-1 < \alpha \leq n, \beta \leq n-1. \end{cases}$$

**Definition 3.** The Laplace transform of the continuous function $v(t)$ is defined by

$$V(s) = L[v(t)] = \int_0^\infty e^{-st} v(t) dt \text{ , where } s \text{ is a real or complex number.}$$

**Definition 4.** The Laplace transform $L[v(x,t)]$ of the Caputo fractional derivative is defined as

$$L[D_t^{n\alpha} v(x,t)] = s^{n\alpha} L[v(x,t)] - \sum_{k=0}^{n-1} s^{n\alpha-k-1} v^k(x,0), \ n-1 < n\alpha \leq n.$$

## 3. A computational technique

In this section, we present the basic theory and solution procedure of proposed q-HATM for non-linear partial differential equations. We take a general fractional non-linear non-homogeneous differential equation:

$$D_t^\alpha v(x,t) + R v(x,t) + N v(x,t) = g(x,t), \quad n-1 < \alpha \leq n, \tag{4}$$

where $D_t^\alpha v(x,t)$ denotes fractional derivative of the function $v(x,t)$ due to Caputo, $R$ specifies the linear differential operator, whereas $N$ designates nonlinear differential operator and $g(x,t)$ is the source term.

First, by applying the Laplace transform operator on both sides of Eq. (4), we have

$$s^\alpha L[v(x,t)] - \sum_{k=0}^{n-1} s^{\alpha-k-1} v^{(k)}(x,0) + L[Rv(x,t)] + L[Nv(x,t)] = L[g(x,t)]. \quad (5)$$

On simplifying, the above equation reduces to

$$L[v(x,t)] - \frac{1}{s^\alpha} \sum_{k=0}^{n-1} s^{\alpha-k-1} v^{(k)}(x,0) + \frac{1}{s^\alpha}\{L[Rv(x,t)] + L[Nv(x,t)] - L[g(x,t)]\} = 0. \quad (6)$$

Define the nonlinear operator as

$$N[\varphi(x,t;q)] = L[\varphi(x,t;q)] - \frac{1}{s^\alpha} \sum_{k=0}^{n-1} s^{\alpha-k-1} \varphi^k(x,t;q)(0^+)$$

$$+ \frac{1}{s^\alpha}\{L[R\varphi(x,t;q)] + L[N\varphi(x,t;q)] - L[g(x,t)]\}, \quad (7)$$

where $q \in \left[0, \frac{1}{n}\right]$ and $\phi(x,t;q)$ is real function of $x, t$ and $q$.

Now, we can construct a homotopy

$$(1-nq)L[\varphi(x,t;q) - v_0(x,t)] = hqH(x,t)N[\varphi(x,t;q)], \quad (8)$$

where $L$ is Laplace operator, $q \in \left[0, \frac{1}{n}\right]$, $n \geq 1$ is the embedding parameter, $H(x,t)$ denotes a non-zero auxiliary function, $h$ is an auxiliary parameter and which is negative in almost all practical situation, $v_0(x,t)$ is initial approximation of $v(x,t)$, and $\varphi(x,t;q)$ is an unknown function. For the embedding parameter $q = 0$ and $q = \frac{1}{n}$, the following result holds:

$$\varphi(x,t;0) = v_0(x,t) \text{ and } \varphi\left(x,t;\frac{1}{n}\right) = v(x,t), \quad (9)$$

respectively. Consequently, as $q$ increase from $0$ to $\frac{1}{n}$, the solution $\phi(x,t;q)$ transform from initial guess $v_0(x,t)$ to the solution $v(x,t)$. Now expanding the function $\varphi(x,t;q)$ in series form by employing Taylor's theorem about $q$, we have

$$\varphi(x,t;q) = v_0(x,t) + \sum_{m=1}^{\infty} v_m(x,t)q^m, \quad (10)$$

where

$$v_m(x,t) = \frac{1}{m!} \frac{\partial^m \varphi(x,t;q)}{\partial q^m}\Big|_{q=0}. \quad (11)$$

If auxiliary parameter $h$, the initial guess $v_0(x,t)$, and the asymptotic parameter $n$ are properly chosen, the series (10) converges at $q = \frac{1}{n}$. Then we have

$$v(x,t) = v_0(x,t) + \sum_{m=1}^{\infty} v_m(x,t)\left(\frac{1}{n}\right)^m. \quad (12)$$

Define the vector

$$\vec{v}_m = \{v_0(x,t), v_1(x,t), \ldots, v_m(x,t)\}. \quad (13)$$

Differentiating the deformation equation (8) $m$-times with respect to $q$ and then dividing to $m!$ and finally setting $q = 0$, we can construct the $m$-$th$ order deformation equation

$$L[v_m(x,t) - k_m v_{m-1}(x,t)] = hH(x,t)\Re_m(\vec{v}_{m-1}).$$ (14)

Finally, by taking inverse Laplace transform on the both sides of foregoing equation, we get

$$v_m(x,t) = k_m v_{m-1}(x,t) + hL^{-1}[H(x,t)\Re_m(\vec{v}_{m-1})],$$ (15)

where

$$\Re_m(\vec{v}_{m-1}) = \frac{1}{(m-1)!}\frac{\partial^{m-1}N[\varphi(x,t;q)]}{\partial q^{m-1}}|_{q=0}$$ (16)

and

$$k_m = \begin{cases} 0, & m \leq 1, \\ n, & m > 1. \end{cases}$$ (17)

We emphasize that in q-HATM we have great freedom to choose the initial guess $v_0(x,t)$, the non-zero auxiliary parameter $h$, and the asymptotic parameter $n$. Due to the existence of the factor $\left(\frac{1}{n}\right)^m$ in the series solution obtained in the Eq. (12), there is the possibility that there is must faster convergence than can be obtained from the standard HATM. It should be mentioned that for $n = 1$ in Eq. (12), the q-HATM reduces to standard HATM. The auxiliary parameter $h$ plays a very important role in controlling the convergence region and also convergence rate to the solution.

## 4. Numerical examples:

To verify and validate of the present work, we give several examples. These examples are chosen from [21, 23, 24]. Also graphs of comparison between exact and approximate solution were plotted for different values of $h$. Furthermore, we plot the graphs of absolute error function and $h$-curves. Now we began with the following examples.

## Example 4.1.

Consider the following time-fractional linear telegraph equation in the absence of the external source term [21]:

$$\frac{\partial^{2\alpha}v}{\partial t^{2\alpha}} + 2\frac{\partial^{\alpha}v}{\partial t^{\alpha}} + v = \frac{\partial^2 v}{\partial x^2}, \quad 0 < \alpha \leq 1, \quad t \geq 0,$$ (18)

subjected to the initial conditions:

$$v(x,0) = e^x, \ v_t(x,0) = -2e^x, \ 0 < x < 1.$$ (19)

The exact solution of Eq. (18) for $\alpha = 1$ is $v(x,t) = e^{x-2t}$.

By applying the Laplace transform on both sides of Eq. (18) and use the conditions specified in Eq. (19), it leads the following result:

$$L[v] - \left(\frac{e^x}{s} - \frac{2e^x}{s^2}\right) + \frac{1}{s^{2\alpha}}L\left[2\frac{\partial^\alpha v}{\partial t^\alpha} + v - \frac{\partial^2 v}{\partial x^2}\right] = 0. \tag{20}$$

Define the nonlinear operator as

$$N[\varphi(x,t;q)] = L[\varphi(x,t;q)] - \left(\frac{e^x}{s} - \frac{2e^x}{s^2}\right) + \frac{1}{s^{2\alpha}}L\left[2\frac{\partial^\alpha \varphi(x,t;q)}{\partial t^\alpha} + \varphi(x,t;q) - \frac{\partial^2 \varphi(x,t;q)}{\partial x^2}\right]. \tag{21}$$

By using the aforesaid procedure of proposed numerical scheme, we can construct the $m$-$th$ order deformation equation for $H(x,t) = 1$ as

$$L[v_m(x,t) - k_m v_{m-1}(x,t)] = h\Re_m(\vec{v}_{m-1}), \tag{22}$$

where

$$\Re_m(\vec{v}_{m-1}) = L[v_{m-1}] - \left(1 - \frac{k_m}{n}\right)\left(\frac{e^x}{s} - \frac{2e^x}{s^2}\right) + \frac{1}{s^{2\alpha}}L\left[2\frac{\partial^\alpha v_{m-1}}{\partial t^\alpha} + v_{m-1} - \frac{\partial^2 v_{m-1}}{\partial x^2}\right]. \tag{23}$$

Now, by taking inverse Laplace transform for Eq. (22) we have

$$v_m(x,t) = k_m v_{m-1}(x,t) + hL^{-1}\{\Re_m(\vec{v}_{m-1})\}. \tag{24}$$

By using the initial approximations $v_0(x,t) = e^x$, $v_t(x,0) = -2e^x$ and iterative formula (24), we obtain

$$v_0 = e^x(1-2t), \quad v_1 = \frac{-4he^x}{\Gamma(\alpha+2)}t^{\alpha+1}, \quad v_2 = \frac{-4h(n+h)e^x}{\Gamma(\alpha+2)}t^{\alpha+1} + \frac{-8h^2e^x}{\Gamma(2\alpha+2)}t^{2\alpha+1},$$

$$v_3 = \frac{-4h(n+h)^2e^x}{\Gamma(\alpha+2)}t^{\alpha+1} + \frac{-16h^2(n+h)e^x}{\Gamma(2\alpha+2)}t^{2\alpha+1} + \frac{-16h^3e^x}{\Gamma(3\alpha+2)}t^{3\alpha+1}, \dots$$

Proceeding in this manner, we can also compute the rest of components of $v_m(x,t)$ for $m \geq 4$ of the q-HATM solution, and the obtained solution is expressed as

$$v(x,t) = v_0(x,t) + \sum_{m=1}^{\infty} v_m(x,t)\left(\frac{1}{n}\right)^m.$$

In standard case, i.e., for $\alpha = 1$, $n = 1$ and $h = -1$, the series solution converges to the exact solution as $m \to \infty$

$$v(x,t) = e^x\left(1 + (-2)t + \frac{(-2)^2}{2!}t^2 + \frac{(-2)^3}{3!}t^3 + \dots + \frac{(-2)^k}{k!}t^k + \dots\right) = e^{x-2t}.$$

We compare the exact and approximate solutions found by the q-HATM for $m = 3$. **Figs. 1(a)** and **1(b)** gives the surface of the solution $v(x,t)$ for the Ex. 4.1 at $h = -1, n = 1$ and $\alpha = 1$. Also, the absolute error function is presented in **Fig. 1(c)**. Numerical solution obtained by present method is almost similar with exact solution. Here, $3$-$rd$ order approximation is used for evaluating the approximate solution and the effectiveness of the proposed scheme can be enhanced by computing more terms in series solution. **Fig. 2** records the nature of numerical solution for different Brownian motions $\alpha = 0.5, \alpha = 0.6, \alpha = 0.7, \alpha = 0.8, \alpha = 0.9,$ and standard motion $\alpha = 1$ at $x = 1.5, h = -1, n = 1$. In **Fig. 3** different values of convergence control parameter $h$ is considered to minimize residual errors. **Figs. 4** and **5** represents the $h$-curves of time-fractional derivative $n = 1$ and 2 respectively, at the distinct order. In $h$-curves, the horizontal line segment represents the convergence range of the q-HATM solution. We can recognize from **Figs**. **4** and **5** that the convergence range increases as the order of fractional derivative increases and range of convergence is mainly depends as the positive values of $n$. Thus, the presence of the term $\left(\frac{1}{n}\right)^m$ gives faster convergence in this method.

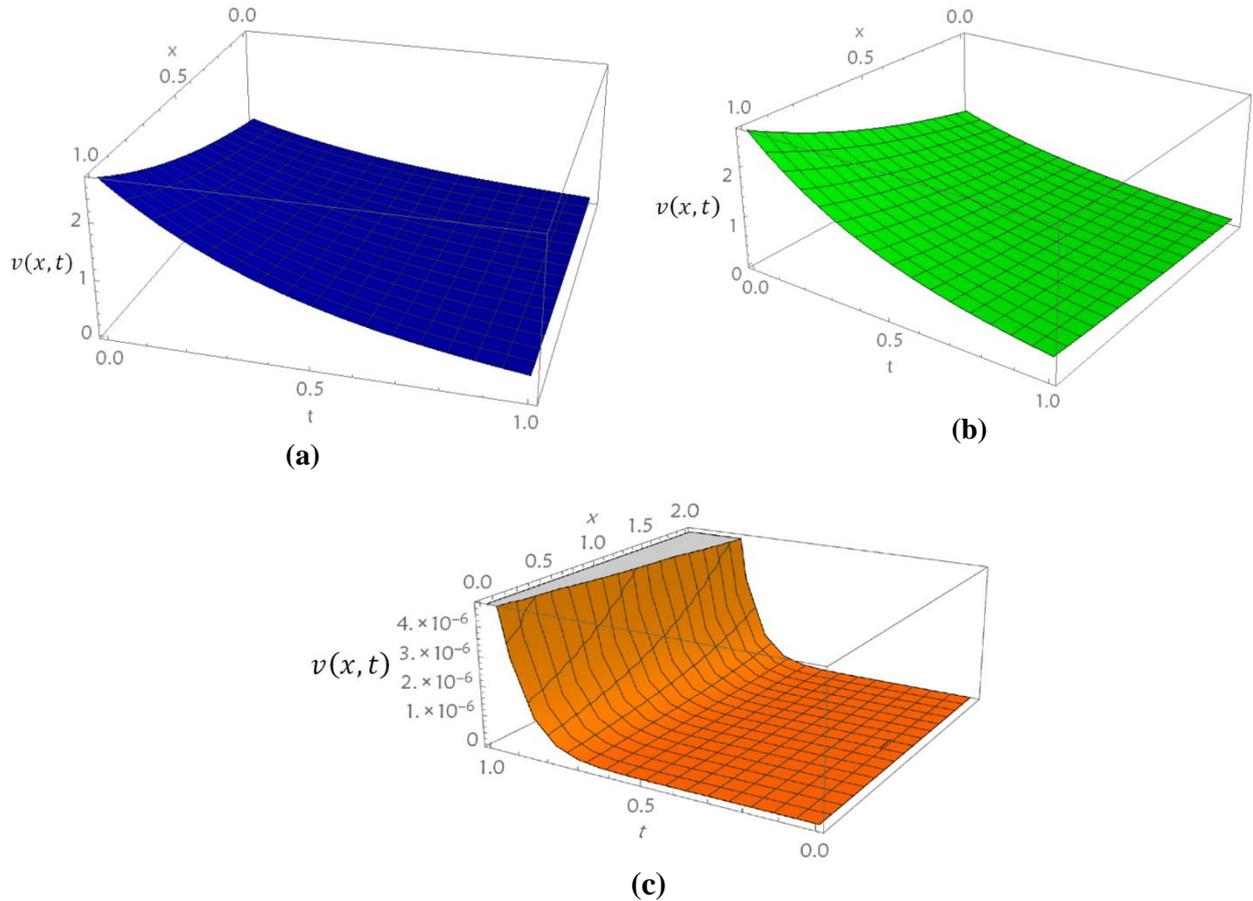

(a)

(b)

(c)

**Fig. 1. (a)** Surface of approximate solution **(b)** Surface of exact solution
**(c)** Surface of absolute error= $\left| v_{exa}\left(x,t\right)-v_{app}\left(x,t\right) \right|$ when $\alpha=1$, $n=1$ and $h=-1$ for Ex. 4.1.

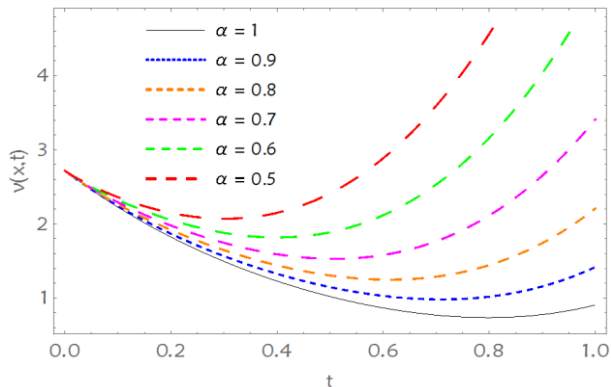

**Fig. 2.** Plot of the q-HATM solution $v(x,t)$ w.r.t $t$ when $h=-1, n=1$ and $x=1.5$ for Ex. 4.1 at different values of $\alpha$.

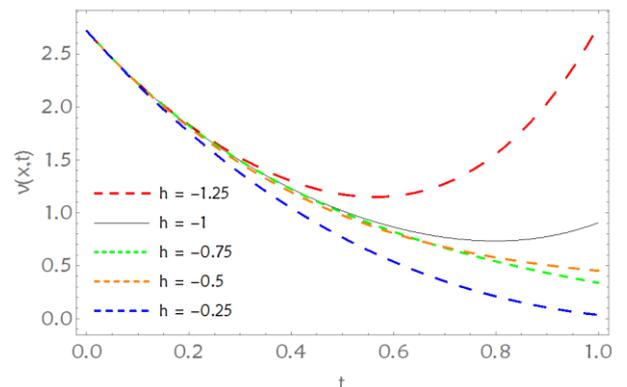

**Fig. 3.** Plot of the q-HATM solution $v(x,t)$ w.r.t $t$ at $\alpha=1$, $n=1$ and $x=1.5$ for Ex. 4.1 at different values of $h$.

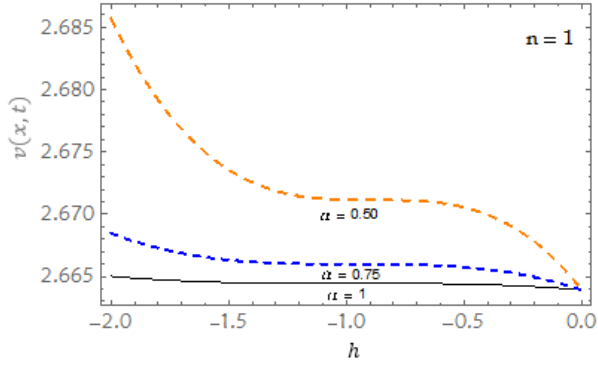
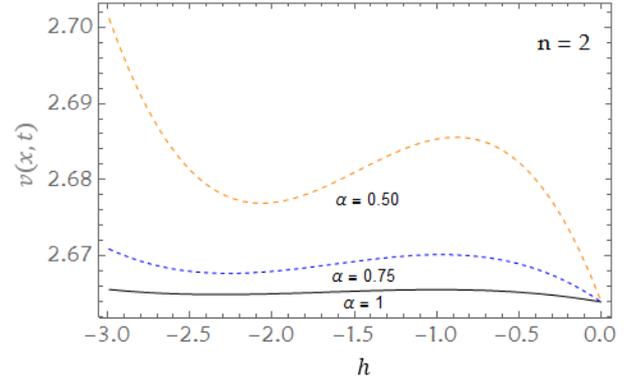

**Fig. 4.** $h$ - curves drwan for the q-HATM solution when $t = 0.01$, $x = 1.5$ and $n = 1$ for Ex. 4.1 at diverse values of $\alpha$.

**Fig. 5.** $h$ - curves drwan for the q-HATM solution when $t = 0.01$, $x = 1.5$ and $n = 2$ for Ex. 4.1 at diverse values of $\alpha$.

### Example 4.2.

Let us consider the following space- fractional telegraph equation [21]:

$$\frac{\partial^\beta v}{\partial x^\beta} = \frac{\partial^2 v}{\partial t^2} + \frac{\partial v}{\partial t} + v, \qquad t \geq 0, \ \ 1 < \beta \leq 2, \tag{25}$$

subjected to the initial and boundary conditions

$$v(x,0) = e^x, \ \ v(0,t) = e^{-t}, \ v_x(0,t) = e^{-t}, \ \ 0 < x < 1. \tag{26}$$

The exact solution of Eq. (25) is $v(x,t) = e^{x-t}$.

Apply the Laplace transformation to Eq. (25) and use the conditions given in Eq. (26) it leads the result

$$L[v] - \frac{e^{-t}}{s} - \frac{e^{-t}}{s^2} - \frac{1}{s^\beta} L\left[\frac{\partial^2 v}{\partial t^2} + \frac{\partial v}{\partial t} + v\right] = 0. \tag{27}$$

Define the nonlinear operator as

$$N[\varphi(x,t;q)] = L[\varphi(x,t;q)] - \left(\frac{e^{-t}}{s} + \frac{e^{-t}}{s^2}\right) - \frac{1}{s^\beta} L\left[\frac{\partial^2 \varphi(x,t;q)}{\partial t^2} + \frac{\partial \varphi(x,t;q)}{\partial t} + \varphi(x,t;q)\right]. \tag{28}$$

By using the aforesaid procedure of analytical technique for $H(x,t) = 1$, we obtain the $m$-$th$ order deformation equation

$$L[v_m(x,t) - k_m v_{m-1}(x,t)] = h\Re_m(\vec{v}_{m-1}), \tag{29}$$

where

$$\Re_m(\vec{v}_{m-1}) = L[v_{m-1}] - \left(1 - \frac{k_m}{n}\right)\left(\frac{e^{-t}}{s} + \frac{e^{-t}}{s^2}\right) - \frac{1}{s^\beta} L\left[\frac{\partial^2 v_{m-1}}{\partial t^2} + \frac{\partial v_{m-1}}{\partial t} + v_{m-1}\right]. \tag{30}$$

Taking inverse Laplace transform on Eq. (29), we get

$$v_m(x,t) = k_m v_{m-1}(x,t) + hL^{-1}\{\Re_m(\vec{v}_{m-1})\}. \tag{31}$$

For solving the foregoing equations, we can construct the iterations

$$v_0 = e^{-t}(1 + x),$$

$$v_1 = -h e^{-t}\left[\frac{x^\beta}{\Gamma(\beta+1)} + \frac{x^{\beta+1}}{\Gamma(\beta+2)}\right],$$

$$v_2 = -h(n+h)e^{-t}\left[\frac{x^\beta}{\Gamma(\beta+1)} + \frac{x^{\beta+1}}{\Gamma(\beta+2)}\right] + h^2 e^{-t}\left[\frac{x^{2\beta}}{\Gamma(2\beta+1)} + \frac{x^{2\beta+1}}{\Gamma(2\beta+2)}\right],$$

$$v_3 = -h(n+h)^2 e^{-t}\left[\frac{x^\beta}{\Gamma(\beta+1)} + \frac{x^{\beta+1}}{\Gamma(\beta+2)}\right] + 2h^2(n+h)e^{-t}\left[\frac{x^{2\beta}}{\Gamma(2\beta+1)} + \frac{x^{2\beta+1}}{\Gamma(2\beta+2)}\right]$$

$$-h^3 e^{-t}\left[\frac{x^{3\beta}}{\Gamma(3\beta+1)} + \frac{x^{3\beta+1}}{\Gamma(3\beta+2)}\right],$$

$$\vdots$$

Working on the same way, the remaining iterations for $v_m(x,t)$ $(m \geq 4)$ can be obtained. Therefore, the series solution of Eq. (25) using q-HATM is given by

$$v(x,t) = v_0(x,t) + \sum_{m=1}^{\infty} v_m(x,t)\left(\frac{1}{n}\right)^m.$$

It can be observed that, by taking $\beta = 2, n = 1$ and $h = -1$ the series solution $\sum_{m=1}^{N} v_m(x,t)\left(\frac{1}{n}\right)^m$ converges to the exact solution

$$v(x,t) = \left(1 + x + \frac{x^2}{2!} + \frac{x^3}{3!} + \cdots + \frac{x^k}{k!} + \cdots\right)e^{-t} = e^{x-t}.$$

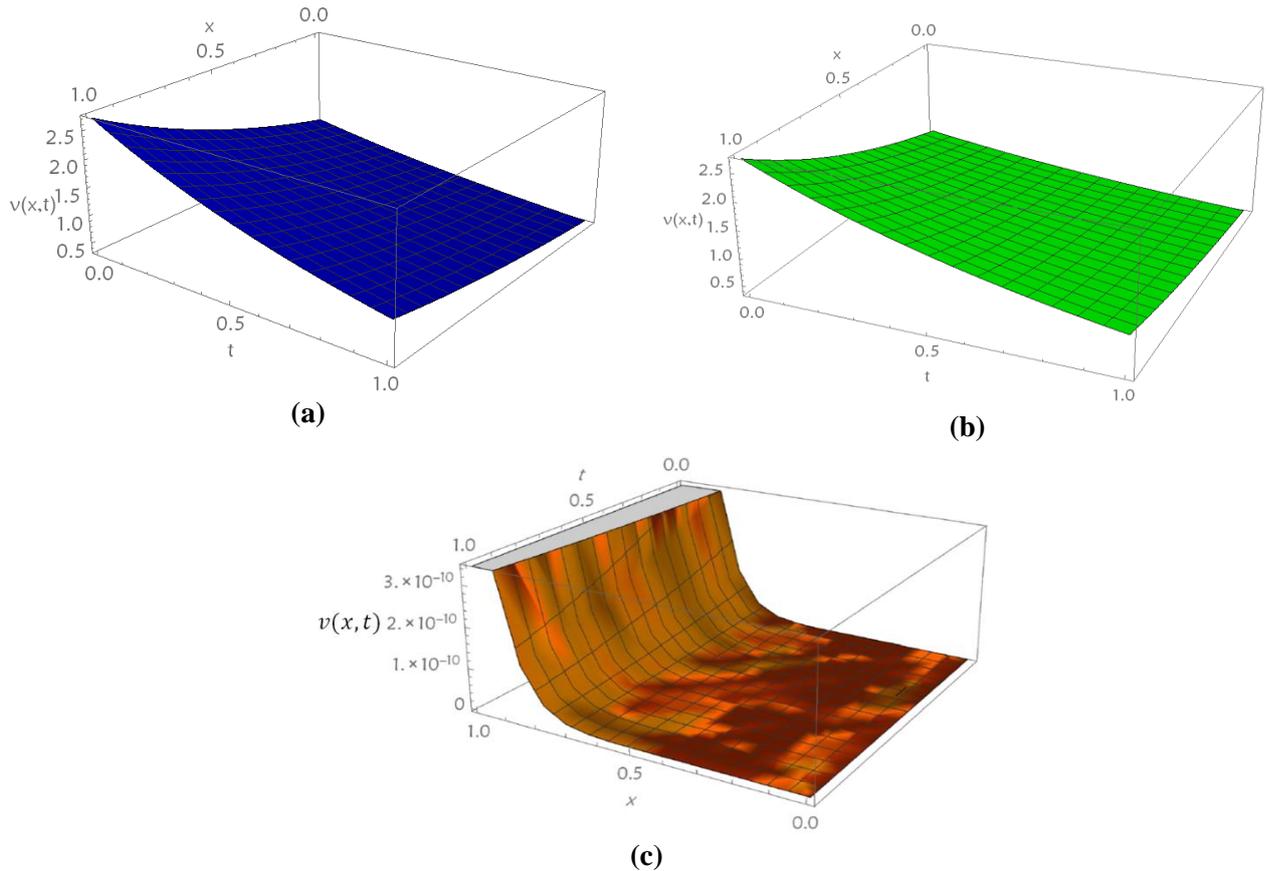

**(a)**　　　　　　　　　　　　　　　**(b)**

**(c)**

**Fig. 6. (a)** Nature of approximate solution, **(b)** Nature of exact solution and **(c)** Nature of Absolute error= $\left|v_{exa.}(x,t) - v_{app.}(x,t)\right|$ when $\beta = 2, n = 1$ and $h = -1$ for Ex. 4.2.

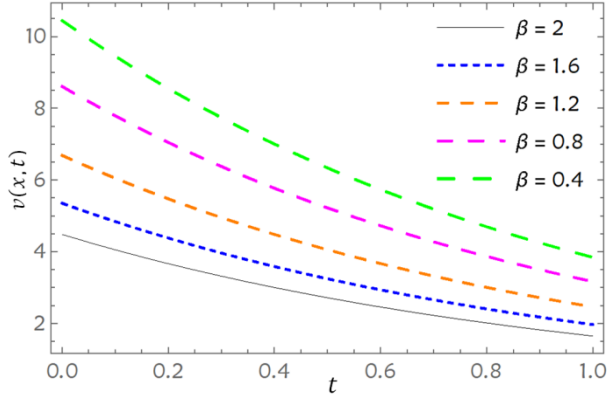

**Fig. 7.** Plot of the q-HATM $v(x,t)$ w.r.t. $t$ at $h = -1, n = 1$ and $x = 1.5$ for Ex. 4.2.

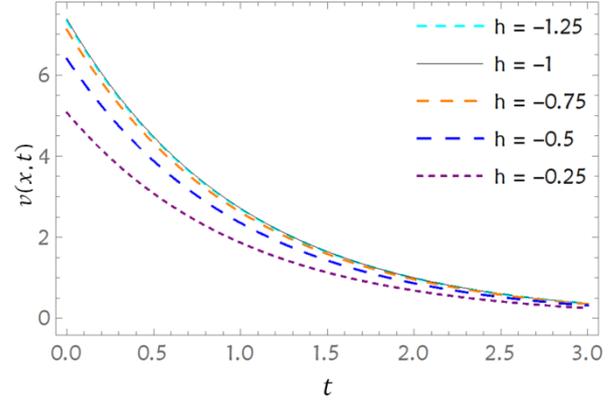

**Fig. 8.** Plot of the q-HATM $v(x,t)$ w.r.t. $t$ at $n = 1$, $\beta = 2$ and $x = 1.5$ for Ex. 4.2.

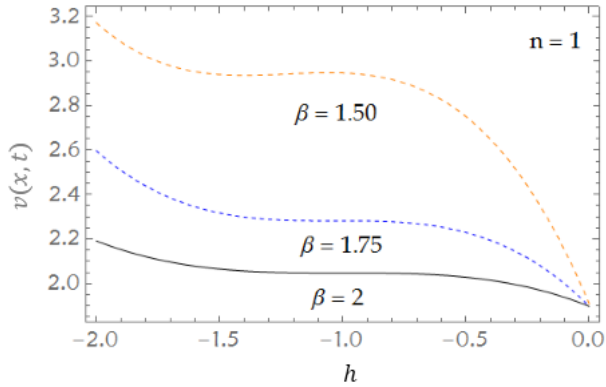

**Fig. 9.** $h$ - curves drwan for the q-HATM solution when $t = 0.01$, $x = 1$ and $n = 1$ for Ex. 4.2 at diverse values of $\beta$.

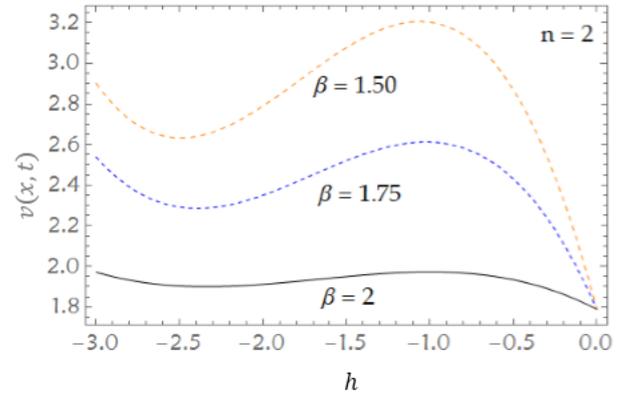

**Fig. 10.** $h$ - curves drwan for the q-HATM solution when $t = 0.01$, $x = 1$ with $n = 2$ for Ex. 4.2 at diverse values of $\beta$.

**Fig. 6** gives plots of the solution $v(x,t)$ for Eqs. (25-26) for $h = -1, n = 1, \beta = 1$. It can be observed from **Figs. 6 (a)** and **(b)** that the numerical solution obtained by q-HATM is almost similar with exact solution. Also graph of absolute error is demonstrated in **Fig. 6 (c)**. The behaviour of numerical solution is recorded for different values of $\beta$ at $h = -1, n = 1$ and $x = 1.5$ in **Fig. 7**. Similarly, in **Fig. 8** different values of convergence control parameter $h$ is considered to minimize residual error. Finally, in **Figs. 9** and **10**, by plotting the $h$-curves we see the convergence region for $n = 1$ and 2.

**Example 4.3.**

Consider the following non-homogeneous space- fractional telegraph equation [24]:

$$\frac{\partial^{2\beta} v}{\partial x^{2\beta}} = \frac{\partial^2 v}{\partial t^2} + \frac{\partial v}{\partial t} + v - x^2 - t + 1, \qquad t \geq 0, \ \ 0 < \beta \leq 1, \qquad (32)$$

subjected to the initial and boundary conditions

$$v(x,0) = x^2, \ \ v(0,t) = t, \ v_x(0,t) = 0, \ \ 0 < x < 1. \qquad (33)$$

Now, applying the aforesaid technique and we define a non-linear operator as

$$N[\varphi(x,t;q)] = L[\varphi(x,t;q)] - \frac{t}{s} + \frac{2}{s^{2\beta+3}} + \frac{(t-1)}{s^{2\beta+2}}$$

$$- \frac{1}{s^{2\beta}} L\left[\frac{\partial^2 \varphi(x,t;q)}{\partial t^2} + \frac{\partial \varphi(x,t;q)}{\partial t} + \varphi(x,t;q)\right]. \qquad (34)$$

Thus, we obtain the $m\text{-}th$ order deformation equation

$$L[v_m(x,t) - k_m v_{m-1}(x,t)] = h \Re_m(\vec{v}_{m-1}), \qquad (35)$$

where

$$\Re_m(\vec{v}_{m-1}) = L[v_{m-1}] - \left(1 - \frac{k_m}{n}\right)\frac{t}{s} + \frac{2}{s^{2\beta+3}} + \frac{(t-1)}{s^{2\beta+2}} - \frac{1}{s^{2\beta}} L\left[\frac{\partial^2 v_{m-1}}{\partial t^2} + \frac{\partial v_{m-1}}{\partial t} + v_{m-1}\right]. \qquad (36)$$

By taking inverse Laplace transform on Eq. (35), we have

$$v_m(x,t) = k_m v_{m-1}(x,t) + h L^{-1}\{\Re_m(\vec{v}_{m-1})\}. \qquad (37)$$

On solving above equations, we get

$$v_0 = t,$$

$$v_1 = -\frac{2hx^{2\beta}}{\Gamma(2\beta+1)} + \frac{2hx^{2\beta+2}}{\Gamma(2\beta+3)},$$

$$v_2 = -\frac{2h(n+h)x^{2\beta}}{\Gamma(2\beta+1)} + \frac{2h(n+h)x^{2\beta+2}}{\Gamma(2\beta+3)} + \frac{2h^2 x^{4\beta}}{\Gamma(4\beta+1)} - \frac{2h^2 x^{4\beta+2}}{\Gamma(4\beta+3)},$$

$$v_3 = -\frac{2h(n+h)^2 x^{2\beta}}{\Gamma(2\beta+1)} + \frac{2h(n+h)^2 x^{2\beta+2}}{\Gamma(2\beta+3)} + \frac{4h^2(n+h)x^{4\beta}}{\Gamma(4\beta+1)} - \frac{4h^2(n+h)x^{4\beta+2}}{\Gamma(4\beta+3)} - \frac{2h^3 x^{6\beta}}{\Gamma(6\beta+1)} + \frac{2h^3 x^{6\beta+2}}{\Gamma(6\beta+3)},$$

$$\vdots$$

Continuing in the same way, we can get the remaining iterates. Therefore, series solution of Eq. (32) by q-HATM is given by

$$v(x,t) = v_0(x,t) + \sum_{m=1}^{\infty} v_m(x,t)\left(\frac{1}{n}\right)^m.$$

Now, for the standard case ($i.e., \beta = 1$) we have $v(x,t) = t + x^2$, which is the exact solution for telegraph equation.

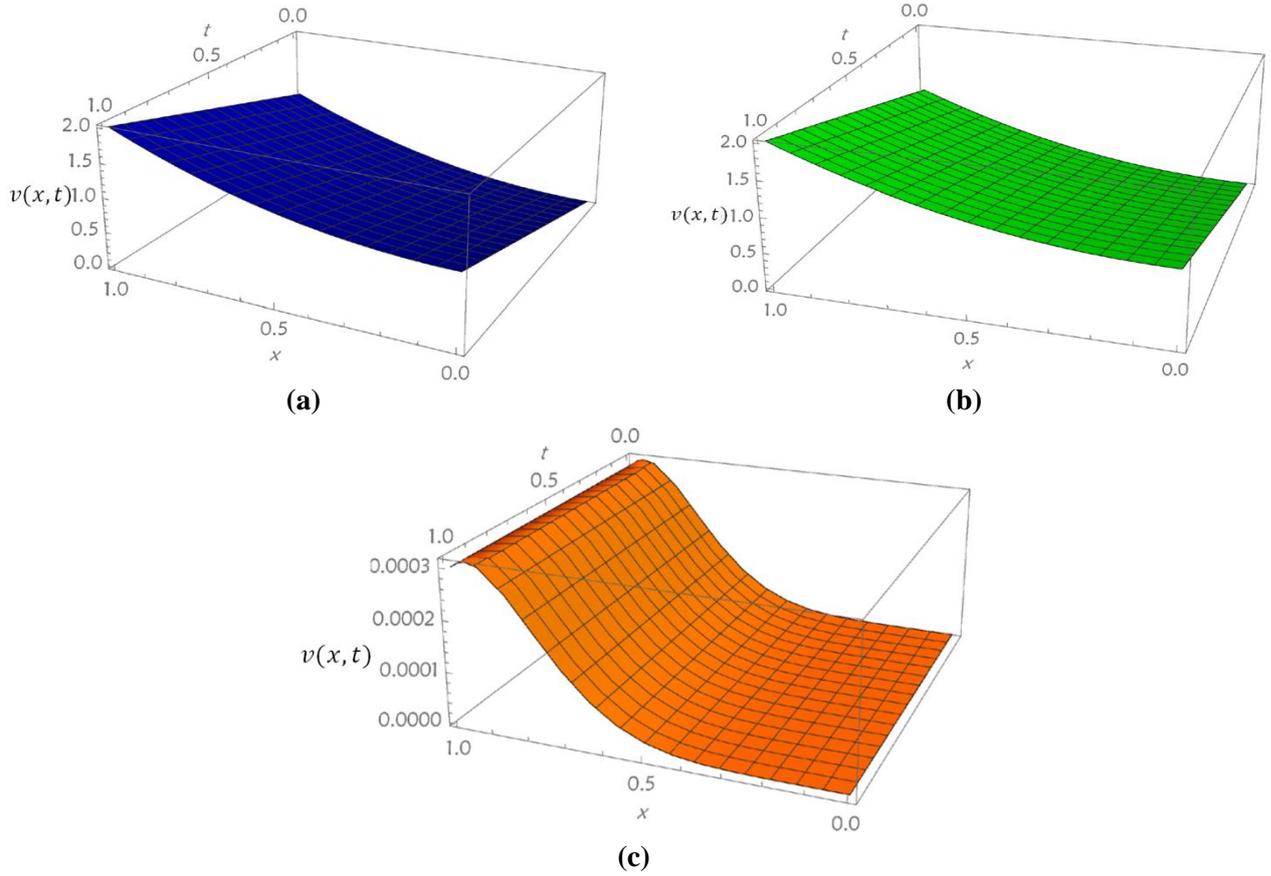

**(a)**　**(b)**

**(c)**

**Fig. 11.** **(a)** Surface of approximate solution, **(b)** Surface of exact solution and **(c)** Surface of absolute error= $\left|v_{exa.}(x,t) - v_{app.}(x,t)\right|$ when $\beta = 1$, $n = 1$ and $h = -1$ for Ex. 4.3.

Comparison between exact and approximate solutions and graph of absolute error function for $m = 3$ are presented in **Figs. 11(a), (b)** and **(c)**, respectively.

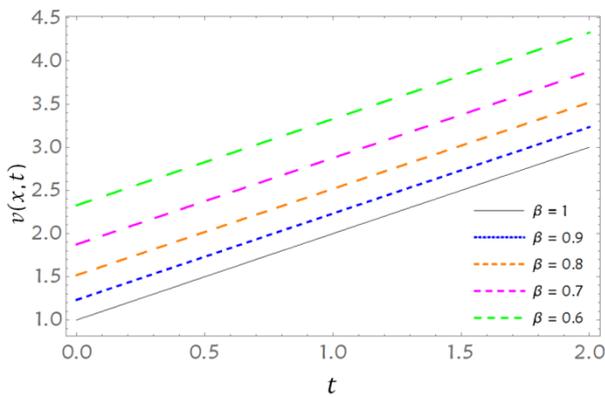

**Fig. 12.** Plot of the q-HATM solution $v(x,t)$ w.r.t. $t$ at $h = -1, n = 1$ and $x = 1.5$ for Ex. 4.3.

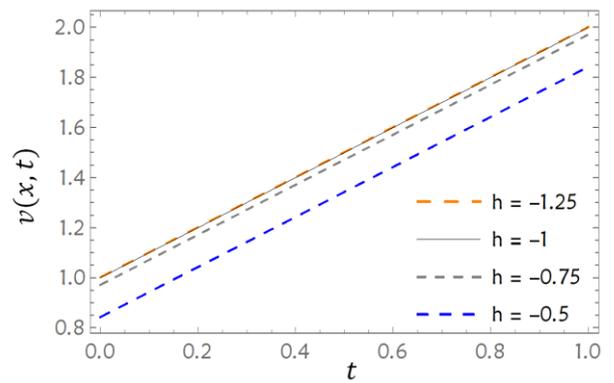

**Fig. 13.** Plot of the q-HATM solution $v(x,t)$ w.r.t. $t$ at $n = 1$, $\beta = 1$ and $x = 2$ for Ex. 4.3.

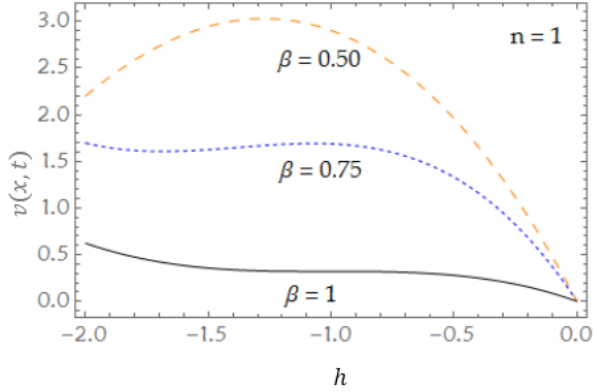
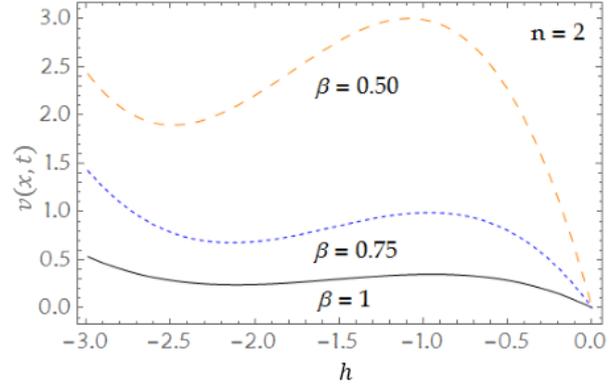

**Fig. 14.** $h$ - curves drwan for the q-HATM solution when $t = 0.01, x = 1$ and $n = 1$ for Ex. 4.3 at different values of $\beta$.

**Fig. 15.** $h$ - curves drwan for the q-HATM solution when $t = 0.01$ and $x = 1$ with $n = 2$ for Ex. 4.3 at different values of $\beta$.

The nature of numerical solution is demonstrated in **Fig. 12** for different values of $\beta$. It can be seen from **Fig. 12** that as the value of $t$ increases q-HATM solution $v(x, t)$ increases. Different values of $h$ is considered to minimize residual error in **Fig. 13**. Also, convergence region for asymptotic parameter $n = 1$ and 2 are shown in **Figs. 14** and **15,** respectively.

**Example 4.4.**

Consider the 2D linear time - fractional Telegraph equation [23]:

$$\frac{\partial^{2\alpha} v}{\partial t^{2\alpha}} + 3 \frac{\partial^{\alpha} v}{\partial t^{\alpha}} + 2v = \frac{\partial^2 v}{\partial x^2} + \frac{\partial^2 v}{\partial y^2}, \quad t \geq 0, \ 0 < \alpha \leq 1, \tag{38}$$

subjected to the initial conditions

$$v(x, y, 0) = e^{x+y}, \quad v_t(x, y, 0) = -3e^{x+y}, \quad 0 < x, y < 1. \tag{39}$$

The exact solution of Eq. (38) is $v(x, y, t) = e^{x+y-3t}$.

By taking Laplace transform on both sides of Eq. (38) and simplifying by the help of Eq. (39), it leads the result

$$L[v] - \left(\frac{e^{x+y}}{s} - \frac{3e^{x+y}}{s^2}\right) + \frac{1}{s^{2\alpha}} L\left[3 \frac{\partial^{\alpha} v}{\partial t^{\alpha}} + 2v - \frac{\partial^2 v}{\partial x^2} - \frac{\partial^2 v}{\partial y^2}\right] = 0. \tag{40}$$

We define the nonlinear operator as

$$N[\varphi(x, y, t; q)] = L[\varphi(x, y, t; q)] - \left(\frac{e^{x+y}}{s} - \frac{3e^x}{s^2}\right)$$

$$+ \frac{1}{s^{2\alpha}} L\left[3 \frac{\partial^{\alpha} \varphi(x,y,t;q)}{\partial t^{\alpha}} + 2\varphi(x, y, t; q) - \frac{\partial^2 \varphi(x,y,t;q)}{\partial x^2} - \frac{\partial^2 \varphi(x,y,t;q)}{\partial y^2}\right]. \tag{41}$$

Using the aforesaid procedure of proposed numerical scheme, the $m$-$th$ order deformation equation for $H(x, y, t) = 1$ is given by

$$L[v_m(x, y, t) - k_m v_{m-1}(x, y, t)] = h\Re_m(\vec{v}_{m-1}), \tag{42}$$

where

$$\Re_m(\vec{v}_{m-1}) = L[v_{m-1}] - \left(1 - \frac{k_m}{n}\right)\left(\frac{e^{x+y}}{s} - \frac{3e^x}{s^2}\right)$$

$$+ \frac{1}{s^{2\alpha}} L\left[3\frac{\partial^\alpha v_{m-1}}{\partial t^\alpha} + 2v_{m-1} - \frac{\partial^2 v_{m-1}}{\partial x^2} - \frac{\partial^2 v_{m-1}}{\partial y^2}\right]. \tag{43}$$

By taking inverse Laplace transform on Eq. (42), we have

$$v_m(x,y,t) = k_m v_{m-1}(x,y,t) + hL^{-1}\{\Re_m(\vec{v}_{m-1})\}. \tag{44}$$

On solving above equations, we get

$$v_0 = e^{x+y}(1-3t), \quad v_1 = \frac{-9he^{x+y}}{\Gamma(\alpha+2)} t^{\alpha+1}, \ v_2 = \frac{9h(n+h)e^{x+y}}{\Gamma(\alpha+2)} t^{\alpha+1} + \frac{-27h^2e^{x+y}}{\Gamma(2\alpha+2)} t^{2\alpha+1},$$

$$v_3 = \frac{9h(n+h)^2e^{x+y}}{\Gamma(\alpha+2)} t^{\alpha+1} + \frac{-54h^2(n+h)e^{x+y}}{\Gamma(2\alpha+2)} t^{2\alpha+1} + \frac{-81h^3e^{x+y}}{\Gamma(3\alpha+2)} t^{3\alpha+1}, \ \ldots$$

Making use of the same process, we can calculate the rest of the components of $v_m(x,t)$. Therefore, series solution of Eq. (38) via q-HATM is given by

$$v(x,y,t) = v_0(x,t) + \sum_{m=1}^{\infty} v_m(x,y,t)\left(\frac{1}{n}\right)^m.$$

It can be observed that, by taking $\alpha = 1, n = 1$ and $h = -1$, the series solution $\sum_{m=1}^{N} v_m(x,y,t)\left(\frac{1}{n}\right)^m$ converges to the exact solution given by

$$u(x,y,t) = e^{x+y}\left(1 + (-3)t + \frac{(-3)^2}{2!}t^2 + \frac{(-3)^3}{3!}t^3 + \frac{(-3)^4}{4!}t^4 + \cdots + \frac{(-3)^k}{k!}t^k + \cdots\right) = e^{x+y-3t}.$$

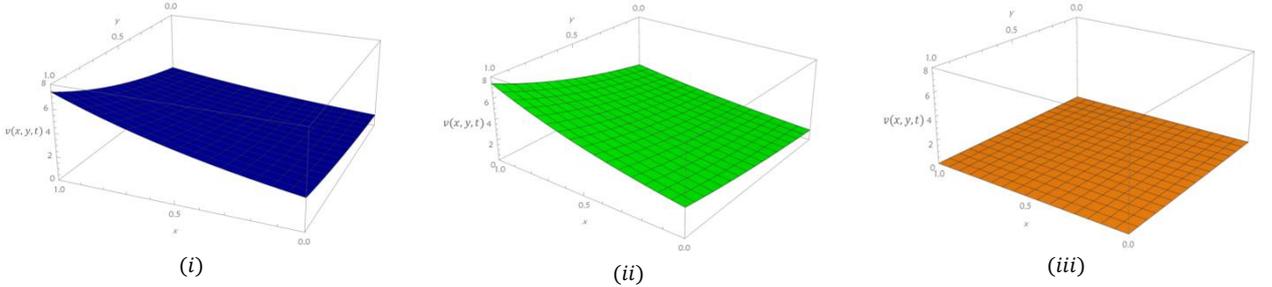

**Fig. 16 (a). (i)** Surface of approximate solution, **(ii)** Surface of exact solution and **(iii)** Surface of Absolute error= $\left|v_{exa.}(x,t) - v_{app.}(x,t)\right|$ when $\alpha = 1, n = 1$ and $h = -1$ for Ex. 4.4 at $\boldsymbol{t = 0.}$

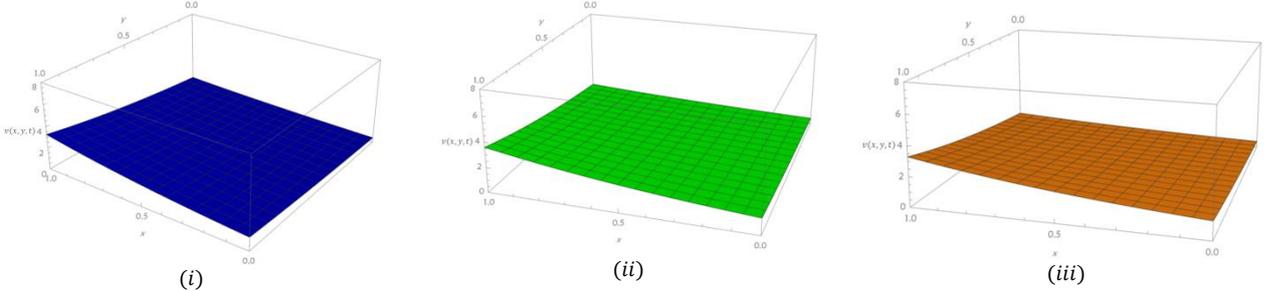

**Fig. 16 (b). (i)** Surface of approximate solution, **(ii)** Surface of exact solution and **(iii)** Surface of Absolute error= $\left|v_{exa.}(x,t) - v_{app.}(x,t)\right|$ when $\alpha = 1, n = 1$ and $h = -1$ for Ex. 4.4 at $\boldsymbol{t = 0.25.}$

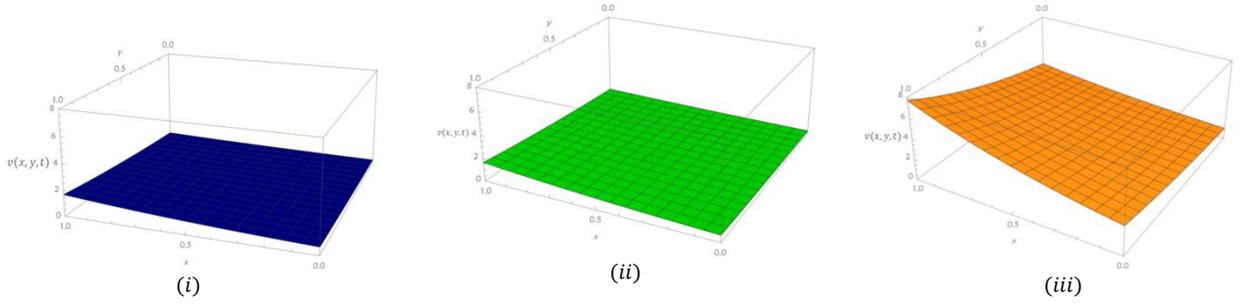

**Fig. 16 (c). (i)** Surface of approximate solution, **(ii)** Surface of exact solution and **(iii)** Surface of Absolute error= $\left| v_{exa.}(x,t) - v_{app.}(x,t) \right|$ when $\alpha = 1, n = 1$ and $h = -1$ for Ex. 4.4 at $\boldsymbol{t = 0.50.}$

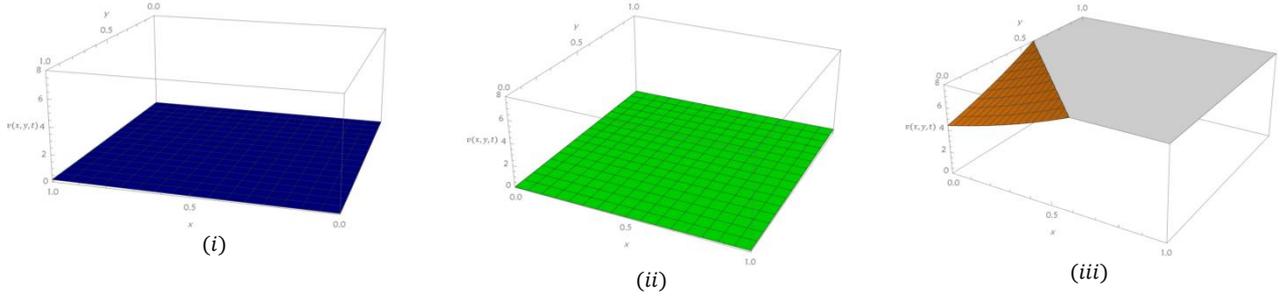

**Fig. 16 (d). (i)** Surface of approximate solution, **(ii)** Surface of exact solution and **(iii)** Surface of Absolute error= $\left| v_{exa.}(x,t) - v_{app.}(x,t) \right|$ when $\alpha = 1, n = 1$ and $h = -1$ for Ex. 4.4 at $\boldsymbol{t = 1.}$

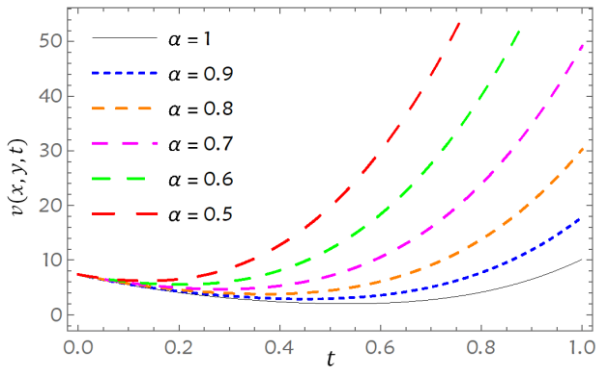

**Fig. 17.** Plot of $v(x,y,t)$ w.r.t. t when $x = 1, y = 1, n = 1$ and $h = -1$ for Ex. 4.4.

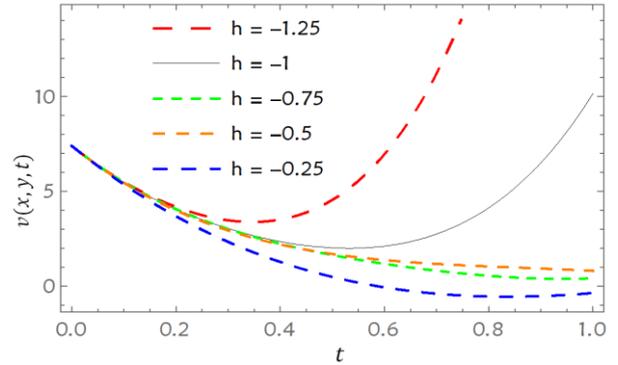

**Fig. 18.** Plot for q-HATM solution $v(x,y,t)$ w.r.t. t when $\alpha = 1, x = 1, y = 1$ and $n = 1$ for Ex. 4.4.

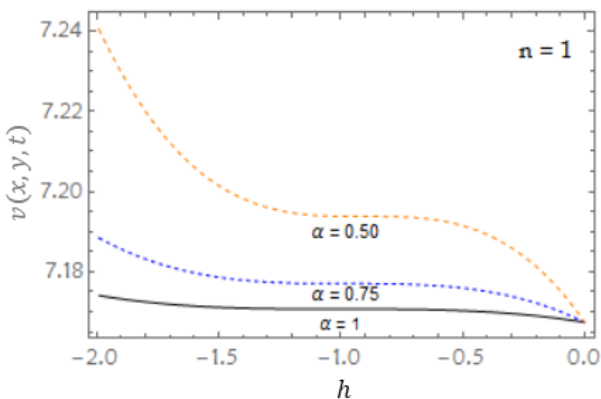

**Fig. 19.** $h$- curves drwan for the q-HATM solution when $x = 1, y = 1, t = 0.01$ and $n = 1$ for Ex. 4.4 with diverse values of $\alpha$.

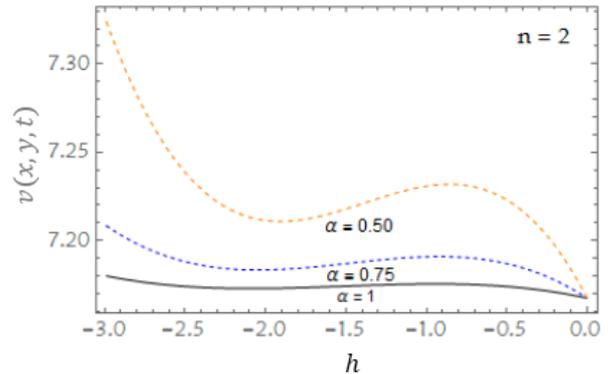

**Fig. 20.** $h$- curves drwan for the q-HATM solution when $x = 1, y = 1, t = 0.01$ and $n = 2$ for Ex. 4.4 with diverse values of $\alpha$.

We plotted the graphs to compare the exact and approximate solution along with absolute error for different intervals of time $t$ in **Figs. 16 (a)** to **16 (d)**. Next, **Fig. 17** exhibits the nature of solution for different values of $\alpha$. In **Fig. 18** different values of convergence control parameter $h$ are picked to minimize residual errors. Finally, $h$-curves are drown for the solution in **Figs. 19** and **20**, respectively for the values of asymptotic parameter $n = 1$ and 2.

**Example 4.5.**

Consider the 3D linear time - fractional Telegraph equation [23]:

$$\frac{\partial^{2\alpha} v}{\partial t^{2\alpha}} + 2\frac{\partial^{\alpha} v}{\partial t^{\alpha}} + 3v = \frac{\partial^2 v}{\partial x^2} + \frac{\partial^2 v}{\partial y^2} + \frac{\partial^2 v}{\partial z^2} \quad t \geq 0, \ 0 < \alpha \leq 1, \tag{45}$$

subjected to the initial conditions

$$v(x, y, z, 0) = sinh(x)\, sinh(y)\, sinh(z),$$

$$v_t(x, y, z, 0) = -sinh(x)\, sinh(y)\, sinh(z), \ 0 < x, y, z < 1. \tag{46}$$

Taking Laplace transform on both sides of Eq. (45) and simplifying, we get

$$L[v] - \left(\frac{1}{s} - \frac{2}{s^2}\right) sinh(x)\, sinh(y)\, sinh(z) + \frac{1}{s^{2\alpha}} L\left[2\frac{\partial^{\alpha} v}{\partial t^{\alpha}} + 3v - \frac{\partial^2 v}{\partial x^2} - \frac{\partial^2 v}{\partial y^2} - \frac{\partial^2 v}{\partial z^2}\right] = 0. \tag{47}$$

Eq. (45) suggests to choose the non-linear operator as

$$N[\varphi(x, y, z, t; q)] = L[\varphi(x, y, z, t; q)] - \left(\frac{1}{s} - \frac{2}{s^2}\right) sinh(x)\, sinh(y)\, sinh(z)$$

$$+ \frac{1}{s^{2\alpha}} L\left[2\frac{\partial^{\alpha} \varphi(x,y,z,t;q)}{\partial t^{\alpha}} + 3\varphi(x, y, t; q) - \frac{\partial^2 \varphi(x,y,z,t;q)}{\partial x^2}\right.$$

$$\left. - \frac{\partial^2 \varphi(x,y,z,t;q)}{\partial y^2} - \frac{\partial^2 \varphi(x,y,z,t;q)}{\partial z^2}\right]. \tag{48}$$

Making use of the stated process of proposed numerical scheme, we obtain the $m^{th}$ order deformation equation for $H(x, y, z, t) = 1$ as

$$L[v_m(x, y, z, t) - k_m v_{m-1}(x, y, z, t)] = h\Re_m(\vec{v}_{m-1}), \tag{49}$$

where

$$\Re_m(\vec{v}_{m-1}) = L[v_{m-1}] - \left(1 - \frac{k_m}{n}\right)\left(\frac{1}{s} - \frac{2}{s^2}\right) sinh(x)\, sinh(y)\, sinh(z)$$

$$+ \frac{1}{s^{2\alpha}} L\left[2\frac{\partial^{\alpha} v_{m-1}}{\partial t^{\alpha}} + 3v_{m-1} - \frac{\partial^2 v_{m-1}}{\partial x^2} - \frac{\partial^2 v_{m-1}}{\partial y^2} - \frac{\partial^2 v_{m-1}}{\partial z^2}\right]. \tag{50}$$

Applying inverse Laplace transform to Eq. (49), it gives

$$v_m(x, y, z, t) = k_m v_{m-1}(x, y, z, t) + hL^{-1}\{\Re_m(\vec{v}_{m-1})\}. \tag{51}$$

On solving above equations, we get

$$v_0 = sinh(x)\, sinh(y)\, sinh(z)\, (1 - 2t),$$

$$v_1 = \frac{-4h\, sinh(x)\, sinh(y)\, sinh(z)}{\Gamma(\alpha+2)}\, t^{\alpha+1},$$

$$v_2 = \frac{-4h(n+h)\sinh(x)\sinh(y)\sinh(z)}{\Gamma(\alpha+2)}\,t^{\alpha+1} + \frac{-8h^2\sinh(x)\sinh(y)\sinh(z)}{\Gamma(2\alpha+2)}\,t^{2\alpha+1},$$

$$v_3 = \frac{-4h(n+h)^2\sinh(x)\sinh(y)\sinh(z)}{\Gamma(\alpha+2)}\,t^{\alpha+1} + \frac{-16h^2(n+h)\sinh(x)\sinh(y)\sinh(z)}{\Gamma(2\alpha+2)}\,t^{2\alpha+1}$$

$$+ \frac{-16h^3\sinh(x)\sinh(y)\sinh(z)}{\Gamma(3\alpha+2)}\,t^{3\alpha+1},$$

$\vdots$

Continuing in the same way, we can get the remaining iterates. Therefore, series solution of Eq. (45) by proposed technique is given by

$$v(x,y,z,t) = v_0(x,y,z,t) + \sum_{m=1}^{\infty} v_m(x,y,z,t)\left(\frac{1}{n}\right)^m.$$

The exact solution of the Eq. (45) is given by

$$u(x,y,z,t) = \sinh(x)\sinh(y)\sinh(z)\left(1 + (-2)t + \frac{(-2)^2}{2!}t^2 + \frac{(-2)^3}{3!}t^3 + \frac{(-2)^4}{4!}t^4 + \cdots\right)$$

$$= e^{-2t}\sinh(x)\sinh(y)\sinh(z).$$

For $\alpha = 1,\ n = 1$ and $h = -1$, the series solution $\sum_{m=1}^{N} v_m(x,y,z,t)\left(\frac{1}{n}\right)^m$ converges to the exact solution.

**Table.** Comparison of numerical results obtained by RDTM and q-HATM with exact solution at different values of $x, y, z$ and $t$ when $\alpha = 1, h = -1$ and $n = 1$ for Ex. 4.5.

| $x = y = z$ | $t$ | $RDTM$ [20] | $q - HATM$ | Exact Solution | Absolute Error |
|---|---|---|---|---|---|
| 0.25 | | 0.0097769 | 0.0097769 | 0.0097772 | 0.0000003 |
| 0.50 | 0.25 | 0.0858202 | 0.0858202 | 0.0858231 | 0.0000029 |
| 0.75 | | 0.3372527 | 0.3372527 | 0.3372641 | 0.0000114 |
| 1 | | 0.9844075 | 0.9844075 | 0.9844404 | 0.0000329 |
| 0.25 | | 0.00591 | 0.00591 | 0.00593 | 0.0000200 |
| 0.50 | 0.50 | 0.05188 | 0.05188 | 0.05205 | 0.0001700 |
| 0.75 | | 0.20388 | 0.20388 | 0.20456 | 0.0006800 |
| 1 | | 0.59512 | 0.59512 | 0.59709 | 0.0019700 |
| 0.25 | | 0.00338 | 0.00338 | 0.00359 | 0.0002100 |
| 0.50 | 0.75 | 0.02973 | 0.02973 | 0.03157 | 0.0018400 |
| 0.75 | | 0.11685 | 0.11685 | 0.12407 | 0.0072200 |
| 1 | | 0.34109 | 0.34109 | 0.36215 | 0.0210600 |
| 0.25 | | 0.00217 | 0.00217 | 0.002181 | 0.0000110 |
| 0.50 | 1 | 0.01912 | 0.01912 | 0.019149 | 0.0000290 |
| 0.75 | | 0.07512 | 0.07512 | 0.07525 | 0.0001300 |
| 1 | | 0.21927 | 0.21927 | 0.21965 | 0.0003800 |

From **Table**, it is observed that $3$-$rd$ order approximation $v_3(x,t)$ is used to compare the numerical results of q-HATM with RDTM and exact solution. An efficiency of q-HATM can be enhanced by computing more approximations.

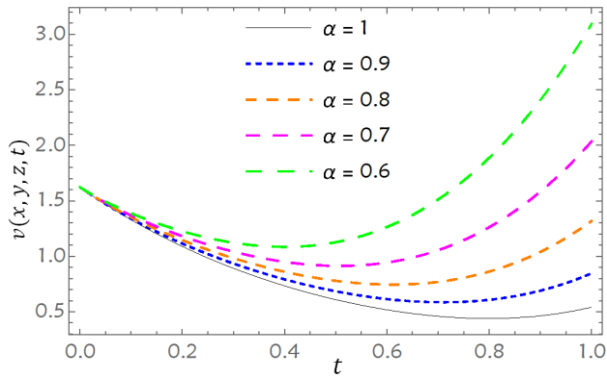

**Fig. 21.** Plot of q-HATM solution $v(x, y, z, t)$ w.r.t. $t$ when $x = 1, y = 1, n = 1, \; z = 1$ and $h = -1$ for Ex. 4.5 for different value of $\alpha$.

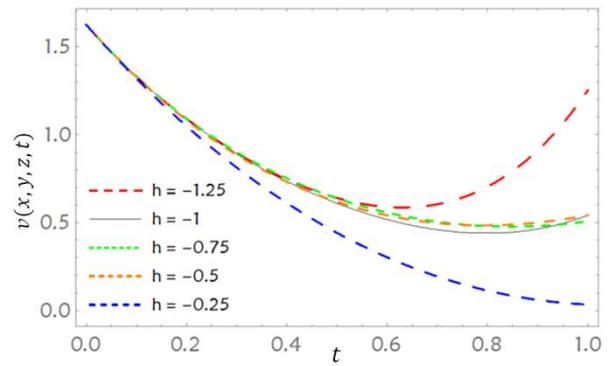

**Fig. 22.** Plot of q-HATM solution $v(x, y, z, t)$ w.r.t. $t$ when $x = 1, y = 1, n = 1, z = 1$ and $\alpha = 1$ for Ex. 4.5 for different value of $h$.

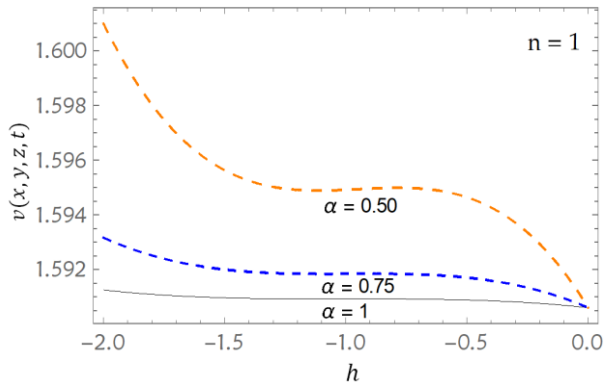

**Fig. 23.** $h$ - curves drwan for the q-HATM solution when $x = 1, y = 1, z = 1, t = 0.01$ and $n = 1$ for Ex. 4.5 with diverse values of $\alpha$.

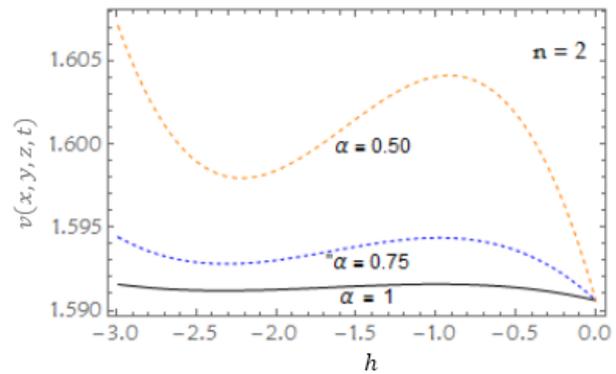

**Fig. 24.** $h$ - curves drwan for the q-HATM solution when $x = 1, y = 1, z = 1, t = 0.01$ and $n = 2$ for Ex. 4.5 with diverse values of $\alpha$.

The behaviour of numerical solution is recorded for different values of $\alpha$ at $h = -1, n = 1$ and $x = y = z = 1$ in **Fig. 21**. Similarly, in **Fig. 22** different values of convergence control parameter $h$ is considered to minimize residual error. Respectively, in **Figs. 23** and **24**, by plotting the $h$-curves we see the convergence region for asymptotic parameter $n = 1$ and 2.

# 6. Conclusion

In this paper, the q-HATM has been successfully applied for getting analytical solutions of multi-dimensional (1D, 2D and 3D) fractional telegraph equations. The main advantage of the method is that it solves the telegraph equation directly without using linearization and perturbation. The q-HATM gives a simple description to adjust and control the convergence of the series solution by selecting the suitable values of the auxiliary parameter $h$ and asymptotic parameter $n$. Also, in this method it is easy to find out the valid regions of $h$ to gain a convergent series solution by mean of so-called $h$-curves. It may be concluded that the q-HATM is easy to implement and very efficient in finding approximate solutions as well as analytical solutions to many fractional physical problems emerging in various fields of science and engineering.

## References

[1] M. Caputo, Elasticita e dissipazione, Zanichelli, Bologna, 1969.

[2] K.B. Oldham, J. Spanier, The fractional calculus, Academic Press, New York, 1974.

[3] K.S. Miller, B. Ross, An introduction to fractional calculus and fractional differential equations, Wiley, New York, 1993.

[4] I. Podlubny, Fractional differential equations, Academic Press, New York, 1999.

[5] A.A. Kilbas, H.M. Srivastava, J.J. Trujillo, Theory and applications of fractional differential equations, Amsterdam, Elsevier, 2006.

[6] A.C. Metaxas, R.J. Meredith, Industrial microwave heating, Peter Peregrinus, London, UK, 1993.

[7] O. Heaviside, Electromagnetic theory, Chelsea Publishing Company, New York, Vol-2 (1899).

[8] V.H. Weston, S. He, Wave splitting of the telegraph equation in $R^3$ and its application to inverse scattering, Inverse Problems, 9 (1993) 789-812.

[9] J. Banasiak, J.R. Mika, Singular perturbed telegraph equations with applications in the random walk theory, J. Appl. Math. Stoch. Anal. 11 (1998) 9-28.

[10] P.M. Jordan, A. Puri, Digital signal propagation in dispersive media, J. Appl. Phys. 85 (1999) 1273-1283.

[11] E. Orsingher, Z. Xuelei, The space fractional telegraph equation and the related fractional telegraph process, Chin. Ann. Math. 24 (2003) 45-56.

[12] S. Momani, Analytic and approximated solutions of space-time fractional telegraph equations, Appl. Math. Compu. 170 (2005) 1126-1134.

[13] R.K. Mohanty, A new unconditionally stable difference schemes for the solution of multi-dimensional telegraphic equations, Int. J. Comput. Math. 86 (2009) 2061-2071.

[14] A. Yildirim, He's homotopy perturbation method for solving the space – and time- fractional telegraph equations, Int. J. Comput. Math. 89 (13) (2010) 2998-3006.

[15] A. Sevimlican, An approximation to solution of space and time fractional telegraph equations by He's variational iteration method, Math. Probl. Eng. (2010), doi:10.1155/2010/290631.

[16] M. Dehghan, S.A. Yousefi, A. Lotfi, The use of He's variational iteration method for solving the telegraph and fractional telegraph equations, Int. J. Numer. Methods Bio. Eng. 27 (2011) 219-231.

[17] R. Jiwari, S. Pandit, R.A. Mittal, A differential quadrature algorithm to solve the two dimensional linear hyperbolic telegraph equation with Dirichlet and Neumann boundary conditions, Appl. Math. Comput. 218 (2012) 7279-7294.

[18] Y. Khan, J. Diblik, N. Faraz, Z. Smarda, An efficient new perturbative Laplace method for space-time fractional telegraph equations, Advances in Difference Equations, 2012, 2012:204.

[19] V.K. Srivastava, M.K. Awasthi, M. Tamsir, RDTM solution of Caputo time fractional-telegraph equation. AIP Adv 2013;3:032142.

[20] F.A. Alawad, E.A. Yousif, A. I. Arbab, A new technique of Laplace variational iteration method for solving space-time fractional telegraph equations, International Journal of Differential Equations, 2013, doi:10.1155/2013/256593.

[21] D. Kumar, J. Singh, S. Kumar, Analytic and approximated solutions of space-time fractional telegraph equations via Laplace transform, Walailak J. Sci. & Tech. 11 (8) (2014) 711-728.

[22] S. Kumar, A new analytical modelling for fractional telegraph equation via Laplace transform, Appl. Math. Mode. 38 (2014) 3154-3163.

[23] V.K. Srivastava, M.K. Awasthi, S. Kumar, Analytical approximations of two and three dimensional time-fractional telegraph equation by reduced differential transform method, Egyptian Journal of Basic and Applied Science, I (2014) 60-66.

[24] A. Prakash, Analytical method for space-fractional telegraph equation by Homotopy perturbation transform method, Nonlinear Engineering, 5 (2) (2016) 123-128.

[25] R.R. Dhunde, G.L. Waghmare, Double Laplace transform method for solving space and time fractional telegraph equations, Int. J. Math. Math. Sci. 2016, doi: 10.1155/2016/1414595.

[26] N. Mollahasani, M.M. Mohseni, K. Afrooz, A new treatment based on hybrid functions to the solution of telegraph equations of fractional order, Appl. Math. Model. 40 (2016) 2804-2814.

[27] R.K. Pandey, H.K. Mishra, Numerical simulation of space-time fractional telegraphs equations with local fractional derivatives via HAFSTM, New Astronomy 57 (2017) 82-93.

[28] A. Prakash, H. Kaur, Numerical solution for fractional model of Fokker-Plank equation by using q-HATM, Chaos, Solitons & Fractals 105 (2017) 99-110.

[29] H.M. Srivastava, D. Kumar, J. Singh, An efficient analytical technique for fractional model of vibration equation, Appl. Math. Model. 45 (2017) 192-204.

[30] D. Kumar, J. Singh, D. Baleanu, A new numerical algorithm for fractional Fitzhugh-Nagumo equation arising in transmission of nerve impulses, Nonlinear Dyn. 91 (2018) 307-317.

[31] J. Singh, D. Kumar, R. Swroop, Numerical solution of time- and space-fractional coupled Burgers' equations via homotopy algorithm, Alexandria Engineering Journal 55 (2) (2016) 1753-1763.